\documentclass[journal]{IEEEtran}
\ifCLASSINFOpdf
\else
\fi
\usepackage{amsmath,amssymb,amstext}
\usepackage{mathtools}
\usepackage[utf8]{inputenc}
\usepackage[T1]{fontenc}
\usepackage{amsfonts}
\newcommand*{\QEDB}{\hfil \ensuremath{\square}}
\usepackage[numbers,sort&compress]{natbib}
\usepackage[hidelinks]{hyperref}

\begin{document}
%
\title{Admissibility and Exact Observability of Observation Operators for Micro-Beam Model: Time and Frequency Domain Approaches}
%
%
%

\author{Mohammad~S.~Edalatzadeh,
        Aria~Alasty,
        and~Ramin~Vatankhah,~\IEEEmembership{Member,~IEEE}
\thanks{Manuscript received August 4, 2015; revised January 15, 2016 and September 26, 2016. This research has been partially supported by Iran's National Science Foundation (INSF).}
\thanks{M.S. Edalatzadeh and A. Alasty (Corresponding Author) are with Department of Mechanical Engineering, Sharif University of Technology, Azadi Ave., Tehran, Iran (e-mails: sajad\_edalat@mech.sharif.edu and aalasti@sharif.edu).}
\thanks{R. Vatankhah, is with the Department of Mechanical Engineering, Shiraz University, Shiraz, Iran (e-mail: rvatankhah@shirazu.ac.ir).}}

%
%

\newcommand{\MYfooter}{\smash{\scriptsize
\hfil\parbox[t][\height][t]{\textwidth}{\centering
0018-9286\,\copyright\,2016 IEEE}\hfil\hbox{}}}

\makeatletter
\def\ps@IEEEtitlepagestyle{%
\def\@oddhead{\mbox{}\scriptsize\rightmark IEEE Transactions on Automatic Control, DOI 10.1109/TAC.2017.2654964, \url{http://ieeexplore.ieee.org/document/7820151/}\hfil \thepage}%
\def\@oddfoot{\MYfooter}}
\makeatother
\pagestyle{IEEEtitlepagestyle}

\makeatletter
\def\ps@headings{%
\def\@oddhead{\mbox{}\scriptsize\rightmark IEEE Transactions on Automatic Control, DOI 10.1109/TAC.2017.2654964, \url{http://ieeexplore.ieee.org/document/7820151/}\hfil \thepage}}
\makeatother
\pagestyle{headings}


\maketitle

\begin{abstract}
This study focuses on the exact observability of a non-classical Euler-Bernoulli micro-beam equation. This non-classical model was derived based on the strain gradient elasticity theory, which is intended to explain the phenomenon of size effect at the micron scale. Spectral properties of the corresponding state operator are studied; an asymptotic expression for eigenvalues is calculated, and eigenfunctions are analyzed in order to check the necessary conditions for the exact observability of the system. By examining the eigenfunctions, it is shown that among non-collocated boundary outputs, only measurement of the non-classical moment at the root of the beam yields an admissible observation operator and also defines an exactly observable system. An alternative proof based on the multiplier method, which is commonly employed in the literature on the observability and controllability of infinite dimensional dynamical systems, is presented to provide a comparison between the time and frequency domain approaches.
\end{abstract}


%
\IEEEpeerreviewmaketitle

\section{Introduction}
%
%
%
%
\IEEEPARstart{O}{ne} of the most important structural components in an atomic force microscope \cite{binnig1986atomic} and many micro-electromechanical systems (MEMS) \cite{pelesko2002modeling} is the micro cantilever beam. High efficiency and simple manufacturing process of micro-cantilever beams make them to play a significant role in MEMS devices.

In the last two decades, many experimental observations in some metals and polymers have demonstrated that the classical continuum mechanics cannot yield accurate static and dynamic models for micro-scale structures \cite{fleck1994strain}. In this way, investigators proposed non-classical continuum theories to accurately predict the static and dynamic behaviors of micro-scale structures. Modified strain gradient theory or briefly strain gradient theory is one of the most successful and inclusive non-classical continuum theories introduced by Lam et al. in 2003 \cite{lam2003experiments}. Recently, this newly established theory has been extensively utilized to study static and dynamic behaviors of the micro-scale beams. Below, some of these works are outlined briefly.

In 2009, Kong et al. \cite{kong2009static} derived a new governing partial differential equation of motion of vibrating Euler-Bernoulli micro-scale beams by using the strain gradient elasticity theory introduced by Lam et al. In a similar way utilized by Wang et al. \cite{wang2010micro} in 2010, a governing equation of motion of strain gradient Timoshenko micro-beams was formulated. Considering the effect of mid-plane stretching, Kahrobaiyan et al. developed a nonlinear Euler-Bernoulli beam model based on strain gradient elasticity theory in 2011 \cite{kahrobaiyan2011nonlinear}. In 2013, forced vibrations of geometrically nonlinear strain gradient Euler-Bernoulli beams were investigated by Vatankhah et al. \cite{vatankhah2013nonlinear} utilizing perturbation techniques. In 2014, by truncating the governing partial differential equation, Edalatzadeh et al. \cite{edalatzadeh2014suppression,edalatzadeh2016boundary} studied stabilization of the non-classical strain gradient Euler-Bernoulli beams, subjected to some nonlinear distributed forces affecting micron and sub-micron structures. Omitting nonlinear terms, but without resorting to model truncation, the authors \cite{vatankhah2013boundary,vatankhah2014exact} investigated the boundary stabilization and exact controllability of the previous beam model. In these studies, well-posedness, stabilizability, and exact controllability of the closed-loop system have been proven by using operator theory and semigroup techniques. In a recent study, Guzmán and Zhu \cite{guzman2015exact} proved the exact controllability of the same micro-beam model considering only one control input.

The concept of observability for infinite dimensional dynamical systems has received considerable attention in recent years. The survey paper of Lagnese \cite{lagnese1991hilbert} and the general exposition of Bensoussan \cite{bensoussan1990general} provided a comprehensive study in this field.  Dolecki and Russell \cite{dolecki1977general} showed that the concept of exact observability for an infinite dimensional dynamical system is dual to that of exact controllability. Furthermore, the duality between an admissible observation operator and an admissible control operator in this framework was introduced by Salamon \cite{salamon1987infinite}. A general necessary condition for exact observability of the systems governed by partial differential equations (PDEs) was then obtained by Russell and Weiss \cite{russell1994general}. To tackle an observability or controllability problem, most researchers adopt a time domain approach in which the governing PDE or its dual counterpart is manipulated in various ways to meet the necessary conditions for the observability or controllability. These ways include the following: multiplier method \cite{komornik1995exact}, microlocal analysis technique \cite{bardos1992sharp}, and nonharmonic Fourier series \cite{avdonin1995families}. The introduction of the Hautus-test for infinite dimensional dynamical systems established a basis for the frequency domain approach, which has rarely been adopted by researcher. Liu et al. \cite{liu2001exponential} obtained a Hautus-type criterion for the exact controllability of systems with bounded input operators; they then applied the criterion to several elastic systems. Recently, a criterion was established for unbounded observation operators with application to the Schrödinger equation \cite{miller2005controllability}. In addition, the Hautus-type test has been further developed to characterize the exact observability of a system only in terms of observation operators and spectral elements of state operators \cite{ramdani2005spectral}.

Although the controllability and control design problems in flexible structure models such as string, rod, beam, and plate models have been addressed in many studies (e.g. see \cite{guzman2015exact,lasiecka1990exact,hansen1995exact,araruna2008controllability,ozer2011exact,ozer2014exact}), the observability and observer design problems of such models have rarely been investigated in the literature. In 2005, backstepping-based infinite dimensional boundary observers were designed by Smyshlyaev and Kirstic \cite{smyshlyaev2005backstepping} to stabilize a class of one-dimensional parabolic PDEs. In 2008, Nguyen \cite{nguyen2008second} designed an observer with a finite number of measurements for second order infinite dimensional systems. Guo et al. studied the exponential stabilization of a one-dimensional wave equation by a boundary controller with collocated and non-collocated observations in \cite{guo2005regularity} and \cite{guo2007stabilization}, respectively. In another study in 2008 \cite{guo2008dynamic}, they proposed a boundary force control and bending strain measurement  to stabilize a classical Euler-Bernoulli beam model. In this paper, the boundary actuator is attached at the free end of the beam while the bending strain observation occurred at the clamped end. In 2011, infinite-dimensional Luenberger-like observers were suggested for a vibrating rotating classical Euler-Bernoulli beam system with constant angular velocity by Li and Xu \cite{li2011infinite}. In 2012, Dogan and Morgul \cite{dogan2011boundary} achieved the same goal as the previous article but by using a backstepping boundary observer.

This paper extends our recent studies on the boundary stabilization and controllability of the non-classical Euler-Bernoulli micro-beam under collocated controls \cite{vatankhah2013boundary,vatankhah2014exact} to address the non-collocated observation problems. Although the controllability and observability analyses of a collocated closed-loop system are performed straightforwardly, the performance of such systems may not be satisfactory \cite{chodavarapu1996noncollocated}. We approached the stabilization problem in \cite{vatankhah2013boundary} by designing a boundary control law and constructing a suitable Lyapunov function for the collocated control system; however, this method is not effective for stabilizing a non-collocated control system. Operator semigroup theory provides some methods such as Riesz basis approach and spectral analysis to deal with such systems.

The rest of this paper is organized as follows. Section 2 provides a brief description of the non-classical Euler-Bernoulli micro-beam model. Section 3 is devoted to the orthonormal basis and spectrum of the beam state operator; an asymptotic expression for eigenvalues is derived in this section for use in the next section. Section 4 discusses the exact observability of the system for various output operators; in addition, it provides a comparison between time and frequency domain approach used to prove the exact observability of the system. Finally, some concluding remarks are given in section 5.

\section{Micro-beam model}
According to Lam et al. \cite{lam2003experiments}, a more realistic, though more complicated, formulation for flexible micro-scale beams can be derived from modified strain gradient elasticity theory using Euler-Bernoulli beam assumptions and Hamilton's principle. As a result, for a micro-cantilever beam with uniform cross-section $A$ and length $L$, the governing PDE of motion and corresponding boundary conditions (BCs) are derived as follows:
\begin{equation}
{K_1}\frac{{{\partial ^4}\hat w}}{{\partial {{\hat x}^4}}} - {K_2}\frac{{{\partial ^6}\hat w}}{{\partial {{\hat x}^6}}} + \rho A\frac{{{\partial ^2}\hat w}}{{\partial {{\hat t}^2}}} = 0,
\end{equation}
\begin{equation}
\left\{ \begin{array}{l} \displaystyle \hat w\left( {0,\hat t} \right) = \frac{{\partial \hat w}}{{\partial \hat x}}(0,\hat t) = \frac{{{\partial ^2}\hat w}}{{\partial {{\hat x}^2}}}(0,\hat t) = 0,\\[2ex]
\displaystyle {K_2}\frac{{{\partial ^5}\hat w}}{{\partial {{\hat x}^5}}}(L,\hat t) - {K_1}\frac{{{\partial ^3}\hat w}}{{\partial {{\hat x}^3}}}(L,\hat t) = \hat F,\\[2ex]
\displaystyle {K_1}\frac{{{\partial ^2}\hat w}}{{\partial {{\hat x}^2}}}(L,\hat t) - {K_2}\frac{{{\partial ^4}\hat w}}{{\partial {{\hat x}^4}}}(L,\hat t) = {{\hat M}^c},\\[2ex]
\displaystyle {K_2}\frac{{{\partial ^3}\hat w}}{{\partial {{\hat x}^3}}}(L,\hat t) = {{\hat M}^{nc}},\end{array} \right.
\end{equation}
where $\hat{x}$ and $\hat{t}$ denote the spatial and time variables, respectively; $\rho $ is the beam density; $\hat{w}(\hat{x},\hat{t})$ indicates the lateral deflection; $\hat{F},\,\,{{\hat{M}}^{c}}$, and ${{\hat{M}}^{nc}}$ can be considered as control inputs and refer to boundary force, moment, and non-classical moment exerted at the tip of the beam, respectively.  In addition,
\begin{equation}
\begin{array}{l}\displaystyle{K_1} = EI + \mu A\left( {2l_0^2 + \frac{8}{{15}}l_1^2 + l_2^2} \right),\\[2ex]
\displaystyle{K_2} = \mu I\left( {2l_0^2 + \frac{4}{5}l_1^2} \right),\end{array}
\end{equation}
where $I$ is the area moment of inertia of the beam cross-section; $E$ and $\mu $ denote Young and shear modulus, respectively; ${{l}_{0}}$, ${{l}_{1}}$, and ${{l}_{2}}$ are additional material constants associated with higher order stress tensors. It can be observed that setting ${{l}_{0}}$, ${{l}_{1}}$, and ${{l}_{2}}$ to zero leads to the classical Euler-Bernoulli beam model. For the sake of brevity, the following dimensionless variables are introduced:
\begin{equation}
\begin{array}{l}\displaystyle w = \frac{{\hat w}}{L},\,\,\,x = \frac{{\hat x}}{L},\,\,\,t = \sqrt {\frac{{{K_1}}}{{\rho A{L^4}}}} \,\hat t,\,\,\,\zeta  = \frac{{{K_2}}}{{{K_1}{L^2}}},\,\,\\[2ex]
\displaystyle F = \frac{{{L^2}}}{{{K_1}}}\hat F,\,\,{M^c} = \frac{{{L}}}{{{K_1}}}{{\hat M}^c},\,\,{M^{nc}} = \frac{{{1}}}{{{K_1}}}{{\hat M}^{nc}} \cdot \end{array}
\end{equation}
From now on, all variables are dimensionless. By applying the dimensionless variables to the governing equations, the following PDE and corresponding BCs are obtained:
\begin{equation}
\frac{{{\partial ^4}w}}{{\partial {x^4}}} - \zeta \frac{{{\partial ^6}w}}{{\partial {x^6}}} + \frac{{{\partial ^2}w}}{{\partial {t^2}}} = 0.
\end{equation}
\begin{equation}
\left\{ \begin{array}{l}\displaystyle w\left( {0,t} \right) = \frac{{\partial w}}{{\partial x}}\left( {0,t} \right) = \frac{{{\partial ^2}w}}{{\partial {x^2}}}\left( {0,t} \right) = 0,\\[2ex]
\displaystyle \zeta \frac{{{\partial ^5}w}}{{\partial {x^5}}}\left( {1,t} \right) - \frac{{{\partial ^3}w}}{{\partial {x^3}}}\left( {1,t} \right) = F,\,\,\,\,\,\,\,\,\,\,\,\\[2ex]
\displaystyle \frac{{{\partial ^2}w}}{{\partial {x^2}}}\left( {1,t} \right) - \zeta \frac{{{\partial ^4}w}}{{\partial {x^4}}}\left( {1,t} \right) = {M^c},\,\,\,\,\,\,\\[2ex]
\displaystyle \zeta \frac{{{\partial ^3}w}}{{\partial {x^3}}}\left( {1,t} \right) = {M^{nc}}.\,\,\,\,\,\,\,\,\,\,\,\,\,\,\,\,\,\,\,\,\end{array} \right.
\end{equation}
Well-posedness of the above PDE with corresponding BCs has been proved by Vatankhah et al. in 2013 \cite{vatankhah2013boundary}. Moreover, referring to Kong et al. \cite{kong2009static}, the dimensionless kinetic energy $K$ and the strain energy $U$ of this physical system can be determined from
\begin{equation}
\begin{array}{l}\displaystyle K = \frac{1}{2}\int_0^1 {{{\left( {\frac{{\partial w}}{{\partial t}}} \right)}^2}\,dx} ,\\[2ex]
\displaystyle U = \frac{1}{2}\int_0^1 {{{{\left( {\frac{{{\partial ^2}w}}{{\partial {x^2}}}} \right)}^2} + \zeta {{\left( {\frac{{{\partial ^3}w}}{{\partial {x^3}}}} \right)}^2}} dx} .\end{array}
\end{equation}
It should be noted that the sum of kinetic and strain energy (i.e. $E=K+U$) is an invariant of the system with zero inputs.

\section{Orthonormal basis}
In order to be able to use the frequency domain approach, the spectral properties of the state operator corresponding to the beam model are investigated in this section. In the first place, let's consider operator ${{A}_{0}}:D({{A}_{0}})(\subset H)\to H$ defined as follows:
\begin{flalign}
&{A_0}f = f^{(4)} - \zeta f^{(6)},& \notag\\
&D({A_0})= \left\{ {f \in {H^6}(0,1) \cap H_E^3(0,1)|{f^{(2)}}(1) - \zeta {f^{(4)}}(1)} \right. & \notag \\
&\qquad \qquad \; \left. {\,\, = \zeta {f^{(5)}}(1) - {f^{(3)}}(1) = \zeta {f^{(3)}}(1) = 0} \right\},
\end{flalign}
where ${{f}^{(m)}}$ stands for the $m$th order derivative of $f$ with respect to $x$; the Sobolev space ${{H}^{k}}(0,1)$ consists of all functions whose derivatives up to order $k-1$ are absolutely continuous and the $k$th order derivative has finite ${{L}^{2}}$ norm; In addition,
\begin{flalign}
&H \coloneqq {L^2}(0,1) = \{ f:[0,1] \to \mathbb{R} |\, \int_0^1 {{f^2}dx < \infty \} } ,&\\
&H_E^3(0,1) \coloneqq \{ f \in {H^3}(0,1)|\, f(0) = {f^{(1)}}(0) = {f^{(2)}}(0) = 0\} .\notag
\end{flalign}

\textit{Theorem 3.1.} The unbounded operator ${{A}_{0}}$ admits an infinite set of eigenvalues which are positive and increasing; furthermore, the corresponding eigenfunctions form an orthonormal basis of $H$.

The spectral properties of a self-adjoint operator with a compact resolvent are characterized in \cite[Sec. 3.2]{tucsnak2009observation}.  Accordingly, the following lemmas should be proven at first in order to prove the theorem.

\textit{Lemma 3.1.} ${{A}_{0}}$ is a symmetric and strictly positive operator.

\textit{Proof:} For every $f,\,g\in D({{A}_{0}})$, a repeated integration by parts gives
\begin{flalign*}
&\allowdisplaybreaks \left\langle {{A_0}f,g} \right\rangle  = \int_0^1 ({f^{(4)}} - \zeta {f^{(6)}})g\,dx\\
& \allowdisplaybreaks \quad =\left[ {f^{(3)}}g - {f^{(2)}}{g^{(1)}} -\zeta {f^{(5)}}g + \zeta {f^{(4)}}{g^{(1)}} - \zeta {f^{(3)}}{g^{(2)}} 
\right] _0^1\\
& \allowdisplaybreaks \qquad +\int_0^1 {({f^{(2)}}{g^{(2)}} + \zeta {f^{(3)}}{g^{(3)}})\,dx},\\
&\allowdisplaybreaks  \left\langle {f,{A_0}g} \right\rangle = \int_0^1 ({g^{(4)}} - \zeta {g^{(6)}})f\,dx\\
& \allowdisplaybreaks \quad = \left[ {g^{(3)}}f - {g^{(2)}}{f^{(1)}} -\zeta {g^{(5)}}f + \zeta {g^{(4)}}{f^{(1)}} - \zeta {g^{(3)}}{f^{(2)}} \right]_0^1\\
&\allowdisplaybreaks \qquad +\int_0^1 {({g^{(2)}}{f^{(2)}} + \zeta {g^{(3)}}{f^{(3)})}\,dx}. &
\end{flalign*}
According to the definition of the domain of ${{A}_{0}}$, the boundary terms in these equations vanish, giving
\begin{flalign*}
\left\langle {{A_0}f,g} \right\rangle  = \left\langle {f^{(2)},g^{(2)}} \right\rangle  &+ \zeta \left\langle f^{(3)},g^{(3)} \right\rangle  = \left\langle {f,{A_0}g} \right\rangle
\end{flalign*}
which shows that operator ${{A}_{0}}$ is symmetric. It can also be seen that $\left\langle {{A}_{0}}f,f \right\rangle =||{{f}^{(2)}}||_{2}^{2}+\zeta ||{{f}^{(3)}}||_{2}^{2}\,\ge 0$; hence, ${{A}_{0}}$ is indeed a strictly positive operator. \hfill \QEDB

\textit{Lemma 3.2.} ${{A}_{0}}$ is surjective and $A_{0}^{-1}$ is a compact operator.

\textit{Proof:} From the theory of linear ordinary differential equations (ODEs), for any $h\in H$, solving $A_0f=h$ for $f$ leads to
\begin{equation}
f(x) = \int_0^1 {G(x,s)h(s)\,ds} ,
\end{equation}
where $G(x,s)$ is known as the Green’s function, which is the solution of ${{A}_{0}}G(x,s)=\delta (x-s)$. It is derived as:
\[
G(x,s) = \left\{
\begin{aligned} 
&- {c_3}(s) - {c_2}(s) - {\zeta ^{ - 1/2}}\left( {{c_2}(s) - {c_3}(s)} \right){\mkern 1mu} x\\
&\,\,\,\,\, - {\zeta ^{ - 1}}\left( {{c_2}(s) + {c_3}(s)} \right){\mkern 1mu} {x^2}/2 + {c_1}(s){\mkern 1mu} {x^3}\\
&\,\,\,\,\,\,\,\, + {c_2}(s){\mkern 1mu} {e^{x{\zeta ^{ - 1/2}}}} + {c_3}{\mkern 1mu} (s){e^{ - x{\zeta ^{ - 1/2}}}},\quad x \le s,\\
&2{\mkern 1mu} {e^{{\zeta ^{ - 1/2}}}}{\mkern 1mu} \left( {{\zeta ^{ - 1/2}} - 1} \right){c_4}(s){\mkern 1mu}  + {c_6}(s)\\
&\,\,\,\,\, + {c_5}(s)x + {e^{{\zeta ^{ - 1/2}}}}{\mkern 1mu} \left( {{\zeta ^{ - 3/2}} - {\zeta ^{ - 1}}} \right){\mkern 1mu} {c_4}(s){\mkern 1mu} {x^2}\\
&\,\,\,\,\,\,\,\, + {c_4}(s){\mkern 1mu} {\mkern 1mu} {e^{x{\zeta ^{ - 1/2}}}} + {c_4}(s){\mkern 1mu} {e^{(2 - x){\zeta ^{ - 1/2}}}}{\mkern 1mu}, x \ge s.&
\end{aligned} \right.
\]In the above equation, the functions ${{c}_{i}}(s),\,\,i=1,2,..,6$, can be uniquely derived. Briefly, they are adjusted such that at $x=s$ the Green’s function and its derivatives with respect to $x$ up to order four are continuous; in addition, its fifth order derivative must have a jump $-1/\zeta $ at this point. As a result, the Green’s function of the operator ${{A}_{0}}$ has a finite ${{L}^{2}}$ norm; hence, the operator ${{A}_{0}}$ is surjective. Furthermore, $A_{0}^{-1}$ maps $H$ into a dense subset of ${{H}^{6}}(0,1)$ which is compactly embedded in $H$ by the Rellich-Kondrachov compact embedding theorem \cite[Ch. 6]{adams2003sobolev}. Hence, $A_{0}^{-1}$ is compact on $H$, and the proof is complete. \hfill \QEDB

\textit{Proof of Theorem 3.1:} Using Lemmas 3.1, 3.2, and the consequence of the Hilbert-Schmidt theorem for unbounded operators, it can be deduced that ${{A}_{0}}$ generates an infinite set of eigenvalues which are positive and increasing and a set of orthonormal eigenfunctions forming a basis for $H$. \hfill \QEDB

The definition of the operator ${{A}_{0}}$ can help us to express the PDE in (5) and corresponding BCs in (6) (assuming zero inputs) in the form of an evolutionary equation in energy state space $\mathbb{H}=H_E^3(0,1)\times {{L}^{2}}(0,1)$; that is
\begin{equation}
\frac{{d\xi(t)}}{{dt}} = A\xi(t), \quad  \xi(0)=\xi_0\in D(A),
\end{equation}
where $\xi(t)=(w,{{w}_{t}})$, and the operator $A:D(A)(\subset \mathbb{H})\to \mathbb{H}$ is defined as follows:
\begin{equation}
\begin{array}{l}A(f,g) = (g, - {A_0}f),\\D(A) = \left\{ {(f,g)|f \in D({A_0}),\,\,g \in H_E^3(0,1)} \right\}.\end{array}
\end{equation}
The state space $\mathbb{H}$ is equipped with an inner product induced norm defined as
\begin{equation}
\left\| {(f,g)} \right\|_\mathbb{H}^2 = \frac{1}{2}\int_0^1 {\left\{ {{{\left( {{f^{(2)}}} \right)}^2} + \zeta {{\left( {{f^{(3)}}} \right)}^2} + {g^2}} \right\}dx} .
\end{equation}

The operator $A$ is skew-adjoint (i.e. ${{A}^{*}}=-A$) as a result of the symmetry and surjectivity of the operator ${{A}_{0}}$. It can be readily found that the operator $A$ has a compact resolvent due to the compactness of $A_{0}^{-1}$ and the definition of $D(A)$ (see \cite[Proposition 1]{guo2002controllability}). The following proposition gives the relation between eigenvalues and eigenfunctions of the operators ${{A}_{0}}$ and $A$.

\textit{Proposition 3.1.} Consider $\lambda _{i}^{2}$ and ${{\phi }_{i}}$ to be the eigenvalues and corresponding eigenfunctions of the operator ${{A}_{0}}$, respectively. Then, the eigenvalues ${{\mu }_{k}}$ and the corresponding eigenfunctions ${{\psi }_{k}}$ of $A$ for $k\in \mathbb{Z}^*(=\mathbb{Z}\backslash \{0\})$ can be obtained from,
\begin{equation}
\left\{ \begin{array}{l} \displaystyle {\mu _k} = i{\lambda _k}, \quad {\lambda _{ - k}} =  - {\lambda _k},\\[1ex]
\displaystyle {\psi _k} = \frac{1}{{\sqrt 2 }}\left[ {\begin{array}{*{20}{c}}{\frac{1}{{i{\lambda _k}}}{\phi _k}}\\{{\phi _k}}\end{array}} \right],\quad {\phi _{ - k}} =  - {\phi _k}.\end{array} \right.
\end{equation}

\textit{Proof:} Regardless of the definition of the operator ${{A}_{0}}$, a proof of this proposition can be found in \cite[Proposition 3.7.7]{tucsnak2009observation}.\hfill \QEDB

Now, we are able to determine the spectrum of the operator ${{A}_{0}}$ as presented in the following lemma.

\textit{Lemma 3.3.} The spectrum $\sigma ({{A}_{0}})$ of ${{A}_{0}}$ consists of isolated eigenvalues $\lambda_n^2$ which are geometrically simple and for sufficiently large positive integer $n$, admit the following asymptotic expression:
\begin{equation}
{(\frac{{27}}{\zeta }{\lambda_n^2} - \frac{2}{{{\zeta ^3}}})^{\frac{1}{6}}} = \frac{{3\pi }}{{\sqrt 3 }}(n + \frac{1}{2}) + \frac{2}{\pi^2}{ (n + \frac{1}{2})^{ - 2}} + O({n^{ - 3}}). \notag 
\end{equation}

\textit{Proof:} $\lambda^2 \in {{\mathbb{R}}^{+}}$ is an eigenvalue of ${{A}_{0}}$ iff there exist a $\phi \in D({{A}_{0}})$, $\phi \ne 0$ that satisfies
\begin{equation}
\begin{array}{c}{\phi ^{(4)}} - \zeta {\phi ^{(6)}} = {\lambda ^2}\phi ,\\\left\{ \begin{array}{l}\phi (0) = {\phi ^{(1)}}(0) = {\phi ^{(2)}}(0) = 0,\\{\phi ^{(3)}}(1) = 0,\,\,{\phi ^2}(1) - \zeta {\phi ^4}(1) = 0,\\\zeta {\phi ^{(5)}}(1) - {\phi ^{(3)}}(1) = 0.\end{array} \right.\end{array} \label{eigp}
\end{equation}
In order to solve this boundary value problem, roots of the following characteristic polynomial has to be determined:
\begin{equation}
\zeta {s^6} - {s^4} + {\lambda ^2} = 0. \label{charpol}
\end{equation}
which yields six roots $\pm {{s}_{i}},\,i=1,2,3$. Subsequently, a fundamental solution to the ODE (\ref{eigp}) is
\begin{equation}
\phi (x) = {c_1}{e^{{s_1}x}} + {c_2}{e^{ - {s_1}x}} + {c_3}{e^{{s_2}x}} + {c_4}{e^{ - {s_2}x}} + {c_5}{e^{{s_3}x}} + {c_6}{e^{ - {s_3}x}}, \label{eigfun}
\end{equation}
where constants ${{c}_{i}},\,i=1,2,...,6$, can be obtained by applying the BCs in (\ref{eigp}) to $\phi (x)$; doing so leads to the following system of algebraic equations:
\begin{equation}
{\bf{B}}({s_i}){\left[ {{c_1},\,\,{c_2},\,\,{c_3},\,\,{c_4},\,\,{c_5},\,\,{c_6}} \right]^T} = {{\bf{0}}_{6 \times 1}},
\end{equation}
The determinant of $\mathbf{B}({{s}_{i}})$ must be zero so as to have a non-trivial solution for ${{c}_{i}}$. As a result, the eigenvalues can be extracted from the characteristic equation: $\det (\mathbf{B}({{s}_{i}}))=0$. An asymptotic expression of the characteristic equation can be derived for sufficiently large $\lambda $. To this end, the solution ${{s}_{i}}$ to (\ref{charpol}) as $\lambda \to \infty $ are approximated by
\begin{equation}
\begin{array}{l}\displaystyle s_1^2 =  - \frac{1}{{3{\mkern 1mu} \zeta {\mkern 1mu} }}{\left( {27{\zeta ^2}{\lambda ^2} - 2} \right)^{\frac{1}{3}}},\\[2ex]
\displaystyle s_2^2 = \frac{1}{{6{\mkern 1mu} \zeta {\mkern 1mu} }}{\left( {27{\zeta ^2}{\lambda ^2} - 2} \right)^{\frac{1}{3}}} + \frac{{\sqrt 3 {\mkern 1mu} {\mkern 1mu} {\rm{i}}}}{{6\zeta }}{\left( {27{\zeta ^2}{\lambda ^2} - 2} \right)^{\frac{1}{3}}},\\[2ex]
\displaystyle s_3^2 = \frac{1}{{6{\mkern 1mu} \zeta }}{\left( {27{\zeta ^2}{\lambda ^2} - 2} \right)^{\frac{1}{3}}} - \frac{{\sqrt 3 {\mkern 1mu} {\mkern 1mu} {\rm{i}}}}{{6\zeta }}{\left( {27{\zeta ^2}{\lambda ^2} - 2} \right)^{\frac{1}{3}}}.\end{array} \label{si}
\end{equation}
By substituting (\ref{si}) into the characteristic equation and after performing some algebraic manipulations, the asymptotic expression of the characteristic equation for $\zeta =1$ is derived as
\begin{flalign}
&F(\lambda^2 )=-\frac{1}{6}{{a}^{2}}({{e}^{2\bar{q}}}+{{e}^{-2\bar{q}}}+{{e}^{2q}}+{{e}^{-2q}}+8{{e}^{{\bar{q}}}}+8{{e}^{-\bar{q}}}+8{{e}^{q}} \notag \\ 
\label{chareq}
&\, +8{{e}^{-q}}+8{{e}^{\bar{q}-q}}+8{{e}^{q-\bar{q}}}+{{e}^{2\bar{q}-2q}}+{{e}^{2q-2\bar{q}}}+18) \\ 
 &\, \, -m({{e}^{2\bar{q}}}+{{e}^{-2\bar{q}}})-\bar{m}({{e}^{-2q}}+{{e}^{2q}})-2m({{e}^{{\bar{q}}}}+{{e}^{-\bar{q}}}) \notag \\ 
 &\, \, \, -2\bar{m}({{e}^{-q}}+{{e}^{q}})+2({{e}^{\bar{q}-q}}+{{e}^{q-\bar{q}}})+{{e}^{2\bar{q}-2q}}+{{e}^{2q-2\bar{q}}}=0,\notag &
\end{flalign}
where $a={{(27{{\lambda }^{2}}/\zeta -2/{{\zeta }^{3}})}^{1/6}}$, $q=(3/6+\sqrt{3}i/6)a$, and $m=1/2+\sqrt{3}i/2$. In order to find an approximate solution to (\ref{chareq}), all the terms in this equation need to be ordered according to their growth rate; that is
\begin{equation}
\begin{gathered}
\left| {{a}^{2}}({{e}^{2\bar{q}}}+{{e}^{2q}}+8{{e}^{{\bar{q}}}}+8{{e}^{q}}) \right|\in\Theta({{a}^{2}}{{e}^{a}}), \\ 
\left| m{{e}^{2\bar{q}}}+\bar{m}{{e}^{2q}}+2m{{e}^{{\bar{q}}}}+2\bar{m}{{e}^{q}} \right|\in \Theta({{e}^{a}}), \\ 
\vdots
\end{gathered}
\end{equation}
According to Rouché's theorem, the zeroes of the greatest term (i.e. those of order ${{a}^{2}}{{e}^{a}}$) give the exact number of zeroes and also an estimation of zeroes of (\ref{chareq}). Thus, we set
\[{{a}^{2}}({{e}^{2\bar{q}}}+{{e}^{2q}}+8{{e}^{{\bar{q}}}}+8{{e}^{q}})=0,\]
which is simplified to
\begin{equation}
(\frac{\sqrt{3}}{3}a)+8{{e}^{-a/2}}\cos (\frac{\sqrt{3}}{6}a)=0. \label{cos}
\end{equation}
Equation (\ref{cos}) admits the asymptotic solution ${{a}_{n}}=3\pi (n+1/2)/\sqrt{3}+{{\alpha }_{n}}$. In order to estimate ${{\alpha }_{n}}$, the lower order terms (i.e. those of order ${{e}^{a}}$) need to be considered. In other words, by substituting ${{a}_{n}}$ in
\[\begin{aligned}
& -\frac{1}{6}{{a}^{2}}({{e}^{2\bar{q}}}+{{e}^{2q}}+8{{e}^{{\bar{q}}}}+8{{e}^{q}}) \\
&\,\,\,\,\,\,\,\,-m{{e}^{2\bar{q}}}-\bar{m}{{e}^{2q}}-2m{{e}^{{\bar{q}}}}-2\bar{m}{{e}^{q}}=0, \\ 
\end{aligned}\]
which is simplified to
\[\begin{aligned}
  & -\frac{1}{6}{\left(\pi (n+\frac{1}{2})+\frac{\sqrt{3}}{3}{{\alpha }_{n}}\right)^2} \sin (\frac{\sqrt{3}}{3}{{\alpha }_{n}}) \\ 
 & \,\,\,\,\,-\sin (\frac{\sqrt{3}}{3}{{\alpha }_{n}})+\sqrt{3}\cos (\frac{\sqrt{3}}{3}{{\alpha }_{n}})=0, \\ 
\end{aligned}\]
and by using Maclaurin series, the above equation yields the solution ${{\alpha }_{n}}=2{{(n\pi +\pi /2)}^{-2}}+{{\beta }_{n}},\,\,\,{{\beta }_{n}}\in O({{n}^{-3}})$. Therefore, we obtain
\begin{equation}
{{a}_{n}}=\frac{3\pi }{\sqrt{3}}\left(n+\frac{1}{2}\right)+\frac{2}{\pi^2}{\left(n+\frac{1}{2}\right)^{-2}}+{{\beta }_{n}},\; n\to \infty . \label{asym}
\end{equation}
For the purposes of this paper, a better approximation is not needed. 

In order to prove that the eigenvalues are geometrically simple, we first show by contradiction that $\phi^{(3)}(0)$ cannot be zero for an eigenfunction of the operator $A_0$. To this end, multiply the ODE in (\ref{eigp}) by $(x-1){{\phi }^{(1)}}(x)$, integrate over $x$, and perform repeated integration by parts; it then follows that
\begin{flalign}
&\int_{0}^{1}{(x-1){{\phi }^{(1)}}(x)\left[ {{\phi }^{(4)}}(x)-\zeta {{\phi }^{(6)}}(x)-{{\lambda }^{2}}\phi (x) \right]\,\,dx} \notag &\\ 
& \qquad \qquad =\left[ -{{\phi }^{(1)}}(x){{\phi }^{(2)}}(x)+(x-1){{\phi }^{(1)}}(x){{\phi }^{(3)}}(x) \right. \notag &\\
& \qquad \quad-\zeta (x-1){{\phi }^{(1)}}(x){{\phi }^{(5)}}(x) +\zeta {{\phi }^{(1)}}(x){{\phi }^{(4)}}(x) \notag &\\
& \qquad \left. +2\zeta (x-1){{\phi }^{(2)}}(x){{\phi }^{(4)}}(x)-2\zeta {{\phi }^{(2)}}(x){{\phi }^{(3)}}(x) \right]_{0}^{1} \notag &\\ 
&\quad -\frac{1}{2}\left[ (x-1){{\left( {{\phi } ^{(2)}}(x) \right)}^{2}}+\zeta (x-1){{\left( {{\phi }^{(3)}}(x) \right)}^{2}} \right. \notag &\\
&\; \left. +{{\lambda }^{2}}(x-1){{\left( \phi (x) \right)}^{2}} \right]_{0}^{1}+\frac{3}{2}\int_{0}^{1}{{{\left( {{\phi }^{(2)}}(x) \right)}^{2}}}\notag &\\
&+\frac{5\zeta }{2}\int_{0}^{1}{{{\left( {{\phi }^{(3)}}(x) \right)}^{2}}}+\frac{{{\lambda }^{2}}}{2}\int_{0}^{1}{{{\left( \phi (x) \right)}^{2}}}. & \label{integ}
\end{flalign}
Now apply the BCs in (\ref{eigp}) and assume ${{\phi }^{(3)}}(0)=0$, we obtain
\[\frac{3}{2}\int_{0}^{1}{{{\left( {{\phi }^{(2)}}(x) \right)}^{2}}}+\frac{5\zeta }{2}\int_{0}^{1}{{{\left( {{\phi }^{(3)}}(x) \right)}^{2}}}+\frac{{{\lambda }^{2}}}{2}\int_{0}^{1}{{{\left( \phi (x) \right)}^{2}}}=0.\]
The above equation admits a unique solution $\phi (x)=0$, which cannot be an eigenfunction of the invertible operator $A_0$—this argument is suggested as an open problem in Guzmán and Zhu’s paper \cite{guzman2015exact}. Now, let $\phi_1$ and $\phi_2$ be two eigenfunctions of $A_0$ associated with the same eigenvalue $\lambda^2$. Then, the function $\phi(x)=\phi_1^{(3)}(0)\phi_2(x)-\phi_2^{(3)}(0)\phi_1(x)$ satisfies (\ref{eigp}) along with $\phi^{(3)}(0)=0$. As shown, it follows that $\phi(x)=0$, and thus the eigenfunctions $\phi_1$ and $\phi_2$ are not linearly independent. Therefore, the eigenvalues are all geometrically simple, and the proof is complete. \hfill \QEDB
\section{Admissible observation operator and exactly observable system}
This section identifies those observation operators $C:D(A)(\subset \mathbb{H})\to Y$ that are admissible and that define an exactly observable system. Roughly speaking, a system is said to be observable if all states can be determined through some partial measurements of states over a sufficiently long time interval; in addition, the concept of admissibility appears mainly in infinite dimensional dynamical system theory and shows that there exists an output function in ${{L}^{2}}([0,\infty ),Y)$ for any initial state in $\mathbb{H}$. More precisely, the following definition has been introduced.

\textit{Definition 4.1.} The operator $C:D(A)(\subset \mathbb{H})\to Y$ is an admissible observation operator for the semigroup generated by $A$ if there exist two positive constants ${{M}_{1}}$ and $\tau $ such that
\begin{equation}
\int_{0}^{\tau }{\left\| CT(t){{\xi}_{0}} \right\|}_{Y}^{2}dt\le {{M}_{1}}\left\| {{\xi}_{0}} \right\|_{\mathbb{H}}^{2}.
\end{equation}
In addition, the pair $(A,C)$ is exactly observable in time $t\ge \tau $ if there exists positive constant $M_2$ such that
\begin{equation}
\int_{0}^{\tau }{\left\| CT(t){{\xi}_{0}} \right\|_{Y}^{2}}\,dt\ge {{M}_{2}}\left\| {{\xi}_{0}} \right\|_{\mathbb{H}}^{2}.
\end{equation}

Various observation operators can be defined for the system; however, an in-domain point observation may not result in an exactly observable system. Focusing on non-collocated boundary observations, the physical properties that can be measured at the root of the beam are force, moment and non-classical moment. Hence, similar to the governing BCs defined in (6), the observation operators can be defined as follows:
\begin{equation}
\left\{ \begin{aligned}
  & {{C}_{1}}\xi(t)\coloneqq F_{0}^{{}}=\zeta \frac{{{\partial}^{5}}w}{\partial{{x}^{5}}}(0,t)-\frac{{{\partial}^{3}}w}{\partial{{x}^{3}}}(0,t), \\ 
 & {{C}_{2}}\xi(t)\coloneqq M_{0}^{c}=\frac{{{\partial}^{2}}w}{\partial{{x}^{2}}}(0,t)-\zeta \frac{{{\partial}^{4}}w}{\partial{{x}^{4}}}(0,t), \\ 
 & {{C}_{3}}\xi(t)\coloneqq M_{0}^{nc}=\zeta \frac{{{\partial}^{3}}w}{\partial{{x}^{3}}}(0,t). \\ 
\end{aligned} \right. \label{obsv}
\end{equation}

\subsection{Time domain approach}
In the existing literature, admissibility of an observation operator and exact observability of a system are commonly tested by resorting to the basic definition of admissibility and exact observability. In this way, the governing equation is manipulated by performing some integration by parts and using some well-known inequalities in order to construct the desired inequalities in Definition 4.1 (see e.g. \cite{komornik1995exact}, the multiplier method). In what follows, this method is used to prove the admissibility of the observation operator ${{C}_{3}}$ and the exact observability of pair $(A,{{C}_{3}})$.

\textit{Theorem 4.1.} The operator ${{C}_{3}}$ is an admissible observation operator; moreover, the pair $(A,{{C}_{3}})$ is exactly observable.

\textit{Proof:} To prove the theorem, we first multiply the governing equation in (5) by the term $(1-x)\frac{\partial w}{\partial x}$ and integrate with respect to $x$ and $t$; that is
\begin{equation}
\int_{0}^{T}{\int_{0}^{1}{(1-x){{w}_{1}}({{w}_{4}}-\zeta {{w}_{6}}+\ddot{w})dx}}\,dt=0,
\end{equation}
where ${{w}_{m}}$ denotes the $m$th order derivative with respect to $x$, and $\dot{w}$ stands for the derivative of $w$ with respect to time. By performing repeated integration by parts for each integral terms of the above equation and eliminating double integral terms as far as possible, the following equality is derived:
\begin{flalign*}
  & \allowdisplaybreaks \int_{0}^{T}{\left[ (1-x){{w}_{1}}({{w}_{3}}-\zeta {{w}_{5}}) \right]_{0}^{1}dt}+\int_{0}^{T}{\left[ {{w}_{1}}({{w}_{2}}-\zeta {{w}_{4}}) \right]_{0}^{1}dt} \\ 
 & \allowdisplaybreaks \qquad \; \,-\int_{0}^{T}{\int_{0}^{1}{\left[ \frac{1}{2}{{{\dot{w}}}^{2}}+\frac{3}{2}{{\left( {{w}_{2}} \right)}^{2}}+\frac{5\zeta }{2}{{\left( {{w}_{3}} \right)}^{2}} \right]}\,dx}\,dt\\
 & \allowdisplaybreaks \,\,\,\,\,\,\,\,\,\,\,\,+\frac{1}{2}\int_{0}^{T}{\left[ (1-x)w\ddot{w} \right]_{0}^{1}dt}-\frac{1}{2}\left[ \left[ (1-x){{{\dot{w}}}_{1}}w \right]_{0}^{1} \right]_{0}^{T} \\ 
 & \allowdisplaybreaks \,\,\,\,\,\,\,\,\,\,+\int_{0}^{1}{\left[ (1-x){{w}_{1}}\dot{w} \right]_{0}^{T}dx}-\frac{1}{2}\int_{0}^{T}{\left[ \zeta (x-1){{w}_{2}}{{w}_{4}} \right]_{0}^{1}dt} \\ 
 & \allowdisplaybreaks \,\,\,\,\,\,\,\,+\frac{1}{2}\int_{0}^{T}{\left[ \zeta {{w}_{2}}{{w}_{3}} \right]_{0}^{1}dt}+\frac{1}{2}\int_{0}^{T}{\left[ \zeta (x-1){{\left( {{w}_{3}} \right)}^{2}} \right]_{0}^{1}dt} \\ 
 & \allowdisplaybreaks \,\,\,\,\,-\frac{1}{2}\int_{0}^{T}{\left[ 2\zeta {{w}_{3}}{{w}_{2}} \right]_{0}^{1}dt}-\frac{1}{2}\int_{0}^{T}{\left[ \zeta (x-1){{w}_{4}}{{w}_{2}} \right]_{0}^{1}dt} \\
 & \allowdisplaybreaks \,\,+\int_{0}^{T}{\left[ \zeta {{w}_{3}}{{w}_{2}} \right]_{0}^{1}dt}+\frac{3\zeta }{2}\int_{0}^{T}{\left[ {{w}_{2}}{{w}_{3}} \right]_{0}^{1}\,dt}\\ 
 &\allowdisplaybreaks +\frac{1}{2}\int_{0}^{T}{\left[ (x-1){{\left( {{w}_{2}} \right)}^{2}} \right]_{0}^{1}\,dt}=0.& 
\end{flalign*}
Applying the BCs in (6) (inputs are set to zero) to the previous equality reduces this equation to
\begin{flalign}
\int_{0}^{T}{\zeta {{\left. {{\left( {{w}_{3}} \right)}^{2}} \right|}_{x=0}}dt}=&\int_{0}^{T}{\int_{0}^{1}{\left[ {{{\dot{w}}}^{2}}+3{{\left( {{w}_{2}} \right)}^{2}}+5\zeta {{\left( {{w}_{3}} \right)}^{2}} \right]}dx}\,dt \notag \\ 
&+2\int_{0}^{1}{\left[ (x-1){{w}_{1}}\dot{w} \right]_{0}^{T}dx}.& \label{intw3}
\end{flalign}
The next step is to find an upper and lower bound for the right hand side expressions. Focusing on the third integral, one can apply triangular inequality and then Young’s inequality to obtain the following inequality: 
\begin{flalign*}
&\left| 2\int_{0}^{1}{\left[ (x-1){{w}_{1}}\dot{w} \right]_{0}^{T}dx} \right| \le \int_{0}^{1}{\left| x-1 \right| \left[ {{\left( {{w}_{1}} \right)}^{2}}+{{{\dot{w}}}^{2}} \right]_{t=0}dx} \\ 
&\qquad \qquad \qquad \qquad +\int_{0}^{1}{\left| x-1 \right| \left[ {{\left( {{w}_{1}} \right)}^{2}}+{{{\dot{w}}}^{2}} \right]_{t=T}^{{}}dx}.&
\end{flalign*}
Since the maximum value of $|x-1|$ in the interval $[0,1]$ is one, this term can be dropped from the above inequality. Afterwards, a one dimensional version of Poincaré inequality \cite[Lemma 2.1]{krstic2008boundary} can be used to get
\begin{flalign}
&\left| 2\int_{0}^{1}{\left[ (x-1){{w}_{1}}\dot{w} \right]_{0}^{T}dx} \right|\le \int_{0}^{1}{\left[ 4{{\left( {{w}_{2}} \right)}^{2}}+{{{\dot{w}}}^{2}} \right]_{t=0}^{{}}dx} \notag \\ 
 &\hspace{3cm}+\int_{0}^{1}{\left[ 4{{\left( {{w}_{2}} \right)}^{2}}+{{{\dot{w}}}^{2}} \right]_{t=T}^{{}}dx}.& \label{ineq1}
\end{flalign}
The integrals including $w_{2}^{2}$  are bounded by the value of the strain energy of the system; this bound can be found by applying Poincaré inequality to the  expression of the strain energy in (7). That is
\begin{flalign}
&\int_{0}^{1}{{{\left( {{w}_{2}} \right)}^{2}}dx}\le 2\left. {{\left( {{w}_{2}} \right)}^{2}} \right|_{x=0}^{{}}+\int_{0}^{1}{4{{\left( {{w}_{3}} \right)}^{2}}dx} \notag\\ 
 & \qquad\Rightarrow \frac{1}{2}\int_{0}^{1}{\left[ \left( 1+\frac{\zeta }{4} \right){{\left( {{w}_{2}} \right)}^{2}} \right]dx} \le U. \label{ineq2} 
\end{flalign}
Consequently, by substituting (\ref{ineq2}) into (\ref{ineq1}) and using the definition of the kinetic energy in (7), it follows that
\begin{flalign}
&\left| 2\int_{0}^{1}{\left[ (x-1){{w}_{1}}\dot{w} \right]_{0}^{T}dx} \right| \notag \\
&\qquad \qquad \le \left[ \frac{32U}{4+\zeta }+2K \right]_{t=0}^{{}}+\left[ \frac{32U}{4+\zeta }+2K \right]_{t=T}^{{}} \notag \\ 
&\qquad \qquad \qquad \le 2\max (2,\frac{32}{4+\zeta })E.& \label{ineq3}
\end{flalign}
Returning to (\ref{intw3}), for the second integral on the right hand side of this equation, an upper and lower bounds can be readily found as follows:
\begin{equation}
2ET\le \int_{0}^{T}{\int_{0}^{1}{\left( {{{\dot{w}}}^{2}}+3{{\left( {{w}_{2}} \right)}^{2}}+5\zeta {{\left( {{w}_{3}} \right)}^{2}} \right)dx\,dt}}\le 10ET. \label{ineq4}
\end{equation}
Finally, combining (\ref{intw3}), (\ref{ineq3}), and (\ref{ineq4}), we obtain:
\begin{flalign*}
2\left( T-\max (2,\frac{32}{4+\zeta }) \right)E & \le \int_{0}^{T}{\zeta {{\left. {{\left( {{w}_{3}} \right)}^{2}} \right|}_{x=0}}dt} \\ 
&\; \le \left( 10T+2\max (2,\frac{32}{4+\zeta }) \right)E.&
\end{flalign*}
According to the definitions of the observation operator ${{C}_{3}}$, the total energy, the state space, and the corresponding inner product induced norm, the previous inequality can be rewritten as:
\begin{flalign*}
  & 2\zeta \left( T-\max (2,\frac{32}{4+\zeta }) \right)\left\| \xi(0) \right\|_{\mathbb{H}}^{2}\le \int_{0}^{T}{{{\left| {{C}_{3}}\xi(t) \right|}^{2}}dt} \\ 
 & \qquad\qquad\qquad \le \zeta \left( 10T+2\max (2,\frac{32}{4+\zeta }) \right)\left\| \xi(0) \right\|_{\mathbb{H}}^{2},& 
\end{flalign*}
Considering Definition 4.1, it is sufficient to choose ${\tau }=T>\max (2,32/(4+\zeta ))$ to complete the proof. \hfill \QEDB

As can be seen, this usual way of proving an observability estimate is rather constructive and cannot easily be used to show that a system with an observation operator is not exactly observable. In addition, unlike the frequency domain approach, the time domain approach will not yield an optimal observability time ${\tau }$. Accordingly, in the following, the frequency domain approach is adopted to show that the observation operator ${{C}_{1}}$ and ${{C}_{2}}$ are not admissible and that the optimal observability time ${\tau}$ for the observation operator ${{C}_{3}}$ is in fact zero.

\subsection{Frequency domain approach}
Another way to tackle the observability problem for a given observation operator is to consider the image of eigenfunctions of the state operator under the observation operator, providing that the state operator is diagonalizable, which is the case in most physical systems \cite{tucsnak2009observation}. The following proposition provides a powerful tool for studying the observability problem of such systems.

\textit{Proposition 4.1.} Assume that the operator $A$ is skew-adjoint and has a compact resolvent, denoting by ${{\psi }_{k}}$ the eigenfunctions and by $i{{\lambda }_{k}}$ the eigenvalues of $A$ that are simple and ordered such that the sequence ${{\lambda }_{k}}$ is strictly increasing. Then, the operator $C:D(A)(\subset \mathbb{H})\to Y$ is an admissible observation operator for the semigroup generated by $A$, and the pair $(A,C)$ is exactly observable in any time $\tau >0$ if $\lim {}_{\left| k \right|\to \infty }({{\lambda }_{k+1}}-{{\lambda }_{k}})=\infty $ and there exist two positive constants ${{\beta }_{1}}$ and ${{\beta }_{2}}$ such that
\begin{equation}
{{\beta }_{1}}\le {{\left\| C{{\psi }_{k}} \right\|}_{Y}}\le {{\beta }_{2}},\,\,\,\,\,\,\,\,\,\forall k\in \mathbb{Z}^*. \label{Cpsi}
\end{equation}

\textit{Proof:} A proof has been presented in [38, Corollary 6.9.6] by utilizing a wave packets concept. \hfill \QEDB

It is worth mentioning that the inequality ${{\left\| C{{\psi }_{k}} \right\|}_{Y}}\le {{\beta }_{2}}$ solely guarantees the admissibility of an observation operator. 

\textit{Remark 4.1.} For a finite dimensional dynamical system, every observation operator is bounded and hence admissible.

To derive an estimate similar to (\ref{Cpsi}) for our system, we need to determine the asymptotic behavior of the coefficients ${{c}_{i,n}},\,i=1,2,..,6$, in (\ref{eigfun}) corresponding to the eigenfunction $\phi_n(x)$. These asymptotic expressions can be computed with the aid of a symbolic computation package such as MATLAB Symbolic Toolbox. It is observed that
\begin{flalign*}
\|\phi_n(x)\|_2 \in \Theta(a_n^{12}e^{a_n}), \, \phi_n^{(3)}(0)\in \Theta(a_n^{15}e^{a_n}), \\
\phi_n^{(4)}(0)\in \Theta(a_n^{16}e^{a_n}), \, \phi_n^{(5)}(0)\in \Theta(a_n^{17}e^{a_n}),
\end{flalign*}
where $a_n$ is given in (\ref{asym}). Subsequently, images of the eigenfunctions of $A$ under the observation operators (\ref{obsv}) can be obtained by using Proposition 3.1. Then, it is seen that
\begin{gather*}
\underset{\left| k \right|\to \infty }{\mathop{\lim }}\,\left| {{C}_{1}}{{\psi }_{k}} \right|=+\infty, \, \underset{\left| k \right|\to \infty }{\mathop{\lim }}\,\left| {{C}_{2}}{{\psi }_{k}} \right|=+\infty, \\
\underset{\left| k \right|\to \infty }{\mathop{\lim }}\,\left| {{C}_{3}}{{\psi }_{k}} \right|=\sqrt{3\zeta}.
\end{gather*}

These results suggest that the observation operators $C_1$ and $C_2$ are not admissible since the images of the eigenfunctions of $A$ are not bounded under these operators. On the other hand, the observation operator $C_3$ is an admissible observation operator and defines an exactly observable system; this statement is proved in the following theorem. 

\textit{Theorem 4.2.} The observation operator ${{C}_{3}}$ is admissible for the semigroup generated by $A$; moreover, the pair $(A,{{C}_{3}})$ is exactly observable in any time $\tau >0$. 

\textit{Proof:} In the previous section, it is shown that the operator $A$ is skew-adjoint and has a compact resolvent. Furthermore, according to Lemma 3.3, for sufficiently large $k$, the sequence ${{\lambda }_{k}}$ is of order exactly $k^3$, and thus $\lim {}_{\left| k \right|\to \infty }({{\lambda }_{k+1}}-{{\lambda }_{k}})=+\infty $ holds regardless of the choice of an observation operator—this property of the eigenvalues of the micro-beam state operator is suggested as an open problem in Guzmán and Zhu’s paper \cite{guzman2015exact}. 

Now, it is sufficient to show that the sequence $|C_3\psi_k|$ is bounded below and above by some positive numbers. To this end, let $\phi_k$'s be normalized eigenfunctions of $A_0$, apply the BCs in (\ref{eigp}) to (\ref{integ}), we obtain
\begin{equation}
\zeta \left(\phi_k ^{(3)}(0)\right)^2 =\lambda_k^2+ \, \int_0^1{5\zeta\left(\phi_k ^{(3)}(x)\right)^2 +3\left(\phi_k ^{(2)}(x)\right)^2dx}. \label{est1}
\end{equation}
The following integration by parts will then help us to find an estimate for the above integral term:
\begin{flalign*}
&\int_0^1 \phi_k(x)\left[ \phi_k^{(4)}(x)-\zeta \phi_k^{(6)}(x)-\lambda_k ^2 \phi_k(x)\right]dx \notag\\
&\quad=-\lambda_k^2+\zeta\int_0^1{\left(\phi_k ^{(3)}(x)\right)^2}dx +\int_0^1{\left(\phi_k ^{(2)}(x)\right)^2dx}\\
&\qquad+\left[\phi_k(x) \left(\phi_k^{(3)}(x)-\zeta \phi_k^{(5)}(x)\right) -\zeta \phi_k^{(2)}(x)\phi_k^{(3)}(x)\right. \notag\\
&\quad \qquad\; \left.+\phi_k^{(1)}(x)\left(\zeta \phi_k^{(4)}(x)-\phi_k^{(2)}(x)\right)\right]_0^1=0. & \label{integ2}
\end{flalign*}
Applying the BCs in (\ref{eigp}) to the above equation yields
\begin{equation}
\int_0^1{\zeta\left(\phi_k ^{(3)}(x)\right)^2+\left(\phi_k ^{(2)}(x)\right)^2dx}=\lambda_k^2. \label{est2}
\end{equation}
Combining (\ref{est1}) and (\ref{est2}), we obtain 
\begin{equation}
4\lambda_k^2 \le \zeta \left(\phi_k ^{(3)}(0)\right)^2 \le 6\lambda_k^2.
\end{equation} 
Consequently, applying Proposition 3.1, it follows that $\sqrt{2\zeta} \le |C_3\psi_k| \le \sqrt{3\zeta}$. Therefore, according to Proposition 4.1, the observation operator $C_3$ is admissible, and the pair $(A,C_3)$ is exactly observable in any time $\tau>0$. \hfill \QEDB

Knowing that the system is exactly controllable and observable, future research will be focused on designing an exponentially stable observer-base controller. Afterwards, the infinite dimensional observer has to be truncated for practical applications. This late-lumping approach, where an infinite dimensional observer is designed to be reduced to a finite dimensional observer, has several advantages over the early-lumping approach. However, it poses some problems that must be addressed in future work.
\section{Conclusion}
The exact observability of a flexible strain gradient micro-beam was studied with an investigation of different observation operators. It was shown that only the measurement of the non-classical moment at the root of the beam yields an admissible observation operator and defines an exactly observable system. This work contributes to the existing literature by considering a more realistic mathematical model for micro-scale flexible beams as well as adopting and comparing two different approaches to tackling the observability problem.


%





\ifCLASSOPTIONcaptionsoff
  \newpage
\fi



\bibliographystyle{myIEEEtran}

\bibliography{Ref}

\begin{thebibliography}{10}
\providecommand{\url}[1]{#1}
\csname url@samestyle\endcsname
\providecommand{\newblock}{\relax}
\providecommand{\bibinfo}[2]{#2}
\providecommand{\BIBentrySTDinterwordspacing}{\spaceskip=0pt\relax}
\providecommand{\BIBentryALTinterwordstretchfactor}{4}
\providecommand{\BIBentryALTinterwordspacing}{\spaceskip=\fontdimen2\font plus
\BIBentryALTinterwordstretchfactor\fontdimen3\font minus
  \fontdimen4\font\relax}
\providecommand{\BIBforeignlanguage}[2]{{%
\expandafter\ifx\csname l@#1\endcsname\relax
\typeout{** WARNING: IEEEtran.bst: No hyphenation pattern has been}%
\typeout{** loaded for the language `#1'. Using the pattern for}%
\typeout{** the default language instead.}%
\else
\language=\csname l@#1\endcsname
\fi
#2}}
\providecommand{\BIBdecl}{\relax}
\BIBdecl

\bibitem{binnig1986atomic}
G.~Binnig, C.~F. Quate, and C.~Gerber, ``{Atomic Force Microscope},''
  \emph{Phys. Rev. Lett.}, vol.~56, no.~9, pp. 930--933, mar 1986.

\bibitem{pelesko2002modeling}
J.~A. Pelesko and D.~H. Bernstein, \emph{{Modeling MEMS and NEMS}}.\hskip 1em
  plus 0.5em minus 0.4em\relax CRC Press, 2002.

\bibitem{fleck1994strain}
N.~A. Fleck, G.~M. Muller, M.~F. Ashby, and J.~W. Hutchinson, ``{Strain
  gradient plasticity: Theory and experiment},'' \emph{Acta Metallurgica Et
  Materialia}, vol.~42, no.~2, pp. 475--487, 1994.

\bibitem{lam2003experiments}
D.~C.~C. Lam, F.~Yang, A.~C.~M. Chong, J.~Wang, and P.~Tong, ``{Experiments and
  theory in strain gradient elasticity},'' \emph{Journal of the Mechanics and
  Physics of Solids}, vol.~51, no.~8, pp. 1477--1508, 2003.

\bibitem{kong2009static}
S.~Kong, S.~Zhou, Z.~Nie, and K.~Wang, ``{Static and dynamic analysis of micro
  beams based on strain gradient elasticity theory},'' \emph{International
  Journal of Engineering Science}, vol.~47, no.~4, pp. 487--498, 2009.

\bibitem{wang2010micro}
B.~Wang, J.~Zhao, and S.~Zhou, ``{A micro scale Timoshenko beam model based on
  strain gradient elasticity theory},'' \emph{European Journal of Mechanics -
  A/Solids}, vol.~29, no.~4, pp. 591--599, 2010.

\bibitem{kahrobaiyan2011nonlinear}
M.~H. Kahrobaiyan, M.~Asghari, M.~Rahaeifard, and M.~T. Ahmadian, ``{A
  nonlinear strain gradient beam formulation},'' \emph{International Journal of
  Engineering Science}, vol.~49, no.~11, pp. 1256--1267, 2011.

\bibitem{vatankhah2013nonlinear}
R.~Vatankhah, M.~H. Kahrobaiyan, A.~Alasty, and M.~T. Ahmadian, ``{Nonlinear
  forced vibration of strain gradient microbeams},'' \emph{Applied Mathematical
  Modelling}, vol.~37, no. 18–19, pp. 8363--8382, 2013.

\bibitem{edalatzadeh2014suppression}
M.~S. Edalatzadeh, R.~Vatankhah, and A.~Alasty, ``{Suppression of dynamic
  pull-in instability in electrostatically actuated strain gradient beams},''
  in \emph{2014 Second RSI/ISM International Conference on Robotics and
  Mechatronics (ICRoM)}, oct 2014, pp. 155--160.

\bibitem{edalatzadeh2016boundary}
M.~S. Edalatzadeh and A.~Alasty, ``{Boundary exponential stabilization of
  non-classical micro/nano beams subjected to nonlinear distributed forces},''
  \emph{Applied Mathematical Modelling}, vol.~40, no.~3, pp. 2223--2241, feb
  2016.

\bibitem{vatankhah2013boundary}
R.~Vatankhah, A.~Najafi, H.~Salarieh, and A.~Alasty, ``{Boundary stabilization
  of non-classical micro-scale beams},'' \emph{Applied Mathematical Modelling},
  vol.~37, no. 20–21, pp. 8709--8724, 2013.

\bibitem{vatankhah2014exact}
R.~Vatankhah, A.~Najafi, H.~Salarieh, and A.~Alasty, ``{Exact boundary
  controllability of vibrating non-classical Euler–Bernoulli micro-scale
  beams},'' \emph{Journal of Mathematical Analysis and Applications}, vol. 418,
  no.~2, pp. 985--997, 2014.

\bibitem{guzman2015exact}
P.~Guzm{\'{a}}n and J.~Zhu, ``{Exact boundary controllability of a microbeam
  model},'' \emph{Journal of Mathematical Analysis and Applications}, vol. 425,
  no.~2, pp. 655--665, 2015.

\bibitem{lagnese1991hilbert}
J.~E. Lagnese, \emph{{The Hilbert uniqueness method: A retrospective}}.\hskip
  1em plus 0.5em minus 0.4em\relax Berlin, Heidelberg: Springer Berlin
  Heidelberg, 1991, pp. 158--181.

\bibitem{bensoussan1990general}
A.~Bensoussan, ``{On the general theory of exact controllability for skew
  symmetric operators},'' \emph{Acta Applicandae Mathematica}, vol.~20, no.~3,
  pp. 197--229, 1990.

\bibitem{dolecki1977general}
S.~Dolecki and D.~L. Russell, ``{A General Theory of Observation and
  Control},'' \emph{SIAM Journal on Control and Optimization}, vol.~15, no.~2,
  pp. 185--220, 1977.

\bibitem{salamon1987infinite}
D.~Salamon, ``{Infinite-dimensional linear systems with unbounded control and
  observation: a functional analytic approach},'' \emph{Transactions of the
  American Mathematical Society}, vol. 300, no.~2, pp. 383--431, 1987.

\bibitem{russell1994general}
D.~L. Russell and G.~Weiss, ``{A General Necessary Condition for Exact
  Observability},'' \emph{SIAM Journal on Control and Optimization}, vol.~32,
  no.~1, pp. 1--23, 1994.

\bibitem{komornik1995exact}
V.~Komornik, \emph{{Exact Controllability and Stabilization: The Multiplier
  Method}}, ser. Wiley-Masson Series Research in Applied Mathematics.\hskip 1em
  plus 0.5em minus 0.4em\relax Wiley, 1995.

\bibitem{bardos1992sharp}
C.~Bardos, G.~Lebeau, and J.~Rauch, ``{Sharp Sufficient Conditions for the
  Observation, Control, and Stabilization of Waves from the Boundary},''
  \emph{SIAM Journal on Control and Optimization}, vol.~30, no.~5, pp.
  1024--1065, 1992.

\bibitem{avdonin1995families}
S.~A. Avdonin and S.~A. Ivanov, \emph{{Families of Exponentials: The Method of
  Moments in Controllability Problems for Distributed Parameter
  Systems}}.\hskip 1em plus 0.5em minus 0.4em\relax Cambridge University Press,
  1995.

\bibitem{liu2001exponential}
K.~Liu, Z.~Liu, and B.~Rao, ``{Exponential Stability of an Abstract
  Nondissipative Linear System},'' \emph{SIAM Journal on Control and
  Optimization}, vol.~40, no.~1, pp. 149--165, 2001.

\bibitem{miller2005controllability}
L.~Miller, ``{Controllability cost of conservative systems: resolvent condition
  and transmutation},'' \emph{Journal of Functional Analysis}, vol. 218, no.~2,
  pp. 425--444, 2005.

\bibitem{ramdani2005spectral}
K.~Ramdani, T.~Takahashi, G.~Tenenbaum, and M.~Tucsnak, ``{A spectral approach
  for the exact observability of infinite-dimensional systems with skew-adjoint
  generator},'' \emph{Journal of Functional Analysis}, vol. 226, no.~1, pp.
  193--229, 2005.

\bibitem{lasiecka1990exact}
I.~Lasiecka and R.~Triggiani, ``{Exact controllability of the Euler-Bernoulli
  equation with boundary controls for displacement and moment},'' \emph{Journal
  of Mathematical Analysis and Applications}, vol. 146, no.~1, pp. 1--33, 1990.

\bibitem{hansen1995exact}
S.~Hansen and E.~Zuazua, ``{Exact Controllability and Stabilization of a
  Vibrating String with an Interior Point Mass},'' \emph{SIAM Journal on
  Control and Optimization}, vol.~33, no.~5, pp. 1357--1391, 1995.

\bibitem{araruna2008controllability}
F.~D. Araruna and E.~Zuazua, ``{Controllability of the Kirchhoff System for
  Beams as a Limit of the Mindlin–Timoshenko System},'' \emph{SIAM Journal on
  Control and Optimization}, vol.~47, no.~4, pp. 1909--1938, 2008.

\bibitem{ozer2011exact}
A.~{\"{O}}. {\"{O}}zer and S.~W. Hansen, ``{Exact controllability of a Rayleigh
  beam with a single boundary control},'' \emph{Mathematics of Control,
  Signals, and Systems}, vol.~23, no.~1, pp. 199--222, 2011.

\bibitem{ozer2014exact}
A.~{\"{O}}. {\"{O}}zer and S.~W. Hansen, ``{Exact Boundary Controllability
  Results for a Multilayer Rao--Nakra Sandwich Beam},'' \emph{SIAM Journal on
  Control and Optimization}, vol.~52, no.~2, pp. 1314--1337, 2014.

\bibitem{smyshlyaev2005backstepping}
A.~Smyshlyaev and M.~Krstic, ``{Backstepping observers for a class of parabolic
  PDEs},'' \emph{Systems {\&} Control Letters}, vol.~54, no.~7, pp. 613--625,
  2005.

\bibitem{nguyen2008second}
T.~D. Nguyen, ``{Second-order observers for second-order distributed parameter
  systems in R2},'' \emph{Systems {\&} Control Letters}, vol.~57, no.~10, pp.
  787--795, 2008.

\bibitem{guo2005regularity}
B.-Z. Guo and X.~Zhang, ``{The Regularity of the Wave Equation with Partial
  Dirichlet Control and Colocated Observation},'' \emph{SIAM Journal on Control
  and Optimization}, vol.~44, no.~5, pp. 1598--1613, 2005.

\bibitem{guo2007stabilization}
B.~Z. Guo and C.~Z. Xu, ``{The Stabilization of a One-Dimensional Wave Equation
  by Boundary Feedback With Noncollocated Observation},'' \emph{IEEE
  Transactions on Automatic Control}, vol.~52, no.~2, pp. 371--377, feb 2007.

\bibitem{guo2008dynamic}
B.-Z. Guo, J.-M. Wang, and K.-Y. Yang, ``{Dynamic stabilization of an
  Euler–Bernoulli beam under boundary control and non-collocated
  observation},'' \emph{Systems {\&} Control Letters}, vol.~57, no.~9, pp.
  740--749, 2008.

\bibitem{li2011infinite}
X.-D. Li and C.-Z. Xu, ``{Infinite-dimensional Luenberger-like observers for a
  rotating body-beam system},'' \emph{Systems {\&} Control Letters}, vol.~60,
  no.~2, pp. 138--145, 2011.

\bibitem{dogan2011boundary}
M.~Dogan and O.~Morgul, ``{Boundary control of a rotating shear beam with
  observer feedback},'' \emph{Journal of Vibration and Control}, 2011.

\bibitem{chodavarapu1996noncollocated}
P.~A. Chodavarapu and M.~W. Spong, ``{On noncollocated control of a single
  flexible link},'' in \emph{Proceedings of IEEE International Conference on
  Robotics and Automation}, vol.~2, apr 1996, pp. 1101--1106 vol.2.

\bibitem{tucsnak2009observation}
M.~Tucsnak and G.~Weiss, \emph{{Observation and Control for Operator
  Semigroups}}, ser. Birkhäuser Advanced Texts Basler Lehrbücher.\hskip 1em
  plus 0.5em minus 0.4em\relax Birkhäuser Basel, 2009.

\bibitem{adams2003sobolev}
R.~A. Adams and J.~J.~F. Fournier, \emph{{Sobolev Spaces}}, ser. Pure and
  Applied Mathematics.\hskip 1em plus 0.5em minus 0.4em\relax Elsevier Science,
  2003.

\bibitem{guo2002controllability}
B.-Z. Guo and Y.-H. Luo, ``{Controllability and stability of a second-order
  hyperbolic system with collocated sensor/actuator},'' \emph{Systems {\&}
  Control Letters}, vol.~46, no.~1, pp. 45--65, 2002.

\bibitem{krstic2008boundary}
M.~Krstic and A.~Smyshlyaev, \emph{{Boundary Control of PDEs: A Course on
  Backstepping Designs}}, ser. Advances in Design and Control.\hskip 1em plus
  0.5em minus 0.4em\relax SIAM, 2008.

\end{thebibliography}
%

%








\end{document}